\documentclass{amsart}
\usepackage{relsize}
\usepackage{amsmath}
\usepackage{fancyhdr}
\usepackage{amsfonts}
\usepackage{amssymb}
\usepackage{amscd}
\usepackage{amsthm}
\usepackage{mathrsfs}
\newtheorem*{thma}{Theorem 1}
\newtheorem*{thmb}{Theorem 2}
\newtheorem{thm}{Theorem}[section]

\newtheorem{cor}[thm]{Corollary}
\newtheorem{lemma}[thm]{Lemma}
\newtheorem{prop}[thm]{Proposition}

\newtheorem{defi}[thm]{Definition}

\newtheorem{rek}[thm]{Remark}

\newcommand{\be}{\begin{equation}}
\newcommand{\ee}{\end{equation}}
\newcommand{\bea}{\begin{eqnarray}}
\newcommand{\eea}{\end{eqnarray}}
\newcommand{\asd}{\vspace{0.05in}\\\noindent}
\newcommand{\ddd}{\displaystyle}
\newcommand{\mapp}{\mathscr L}
\newcommand{\mapt}{\mathscr T}
\newcommand{\dell}{\partial}
\newcommand{\Ind}{\operatorname{Ind}}
\newcommand{\SL}{\rm{SL}}

\newcommand{\Cc}{\mathcal C}
\newcommand{\Bc}{\mathcal BC}
\newcommand{\ore}{\rm{Re}\,}

\newcommand{\R}{\mathbb R}

\newcommand{\HH}{\mathbb H}

\newcommand{\C}{\mathbb C}

\newcommand{\GG}{\Gamma}
\newcommand{\Ga}{\Gamma}
\newcommand{\Gn}{\tilde\Gamma}
\newcommand{\ga}{\gamma}
\newcommand{\mmod}{\rm{mod}\,}
\newcommand{\Ginf}{\Gamma_{\infty}}
\newcommand{\Gpar}{\Gamma_{\text{par}}}
\newcommand{\Ghyp}{\Gamma_{\text{hyp}}}
\newcommand{\Gell}{\Gamma_{\text{ell}}}
\newcommand{\ze}{Z(s;\Ga;\chi)}
\newcommand{\zen}{Z(s;\Gn;\pi)}
\newcommand{\ziu}{z(s;\Ga;\chi)}
\newcommand{\ziun}{z(s;\Gn;\pi)}
\newcommand{\scat}{\phi(s;\Ga;\chi)}
\newcommand{\scatn}{\phi(s;\Gn;\pi)}
\newcommand{\Sn}{\tilde S_{\infty}}
\newcommand{\lam}{\lambda}
\newcommand{\Linf}{\Lambda_{\infty}}
\newcommand{\om}{\omega}
\newcommand{\Om}{\Omega}
\newcommand{\Omn}{\tilde\Om}
\newcommand{\germa}{\mathfrak a}
\newcommand{\si}{\sigma}
\newcommand{\ep}{\epsilon}
\newcommand{\tr}{{\rm tr}\,}

\newcommand{\eje}{{\bf e}_j}
\newcommand{\mattwo}[4]
{\Big[\begin{array}{cc}
#1&#2\\
#3&#4
\end{array}\Big] }

\begin{document}
\title[Artin Formalism for Selberg Zeta functions]{Artin formalism for Selberg zeta functions of co-finite Kleinian groups}
\author{Eliot Brenner}
\address{University of Minnesota, Minneapolis USA, 55455}
\email{brenn240@umn.edu}
\author{Florin Spinu}
\address{Johns Hopkins University, Baltimore USA, 21218}
\email{fspinu@math.jhu.edu}
\begin{abstract}
Let $\Gamma\backslash\mathbb H^3$ be a finite-volume quotient of the upper-half space,
where $\Gamma\subset {\rm SL}(2,\mathbb C)$ is a discrete subgroup.
To a finite dimensional unitary representation $\chi$ of $\Gamma$ one
associates the Selberg zeta function $Z(s;\Gamma;\chi)$.
In this paper we prove the Artin formalism for the Selberg zeta function.
Namely, if $\tilde\Gamma$ is a finite index group extension of $\Gamma$
in ${\rm SL}(2,\mathbb C)$,
and $\pi={\rm Ind}_{\Gamma}^{\tilde\Gamma}\chi$ is the induced
representation, then $Z(s;\Gamma;\chi)=Z(s;\tilde\Gamma;\pi)$.
In the second part of the paper
we prove by a direct method the analogous identity for the scattering function, namely
$\phi(s;\Gamma;\chi)=\phi(s;\tilde\Gamma;\pi)$,
for an appropriate normalization of the Eisenstein series.

\end{abstract}
\maketitle
\section{Introduction}
\thispagestyle{fancy}
\rhead{}
\renewcommand\headrulewidth{0pt}
\renewcommand\footrulewidth{0.4pt}
\lfoot{{\relsize{-2} MR Classification: 11F72; 11M36.
\newline 
Keywords: Artin Formalism, Selberg Zeta function, Kleinian groups, Fuchsian groups
hyperbolic 3-manifolds, scattering matrix, Eisenstein series}}\
\cfoot{} 

The Artin formalism for the Selberg zeta function associated to
finite area surfaces of constant negative curvature
has been proved by Venkov and al. in \cite{Ve} and \cite{Ve-Zo}.
Two arguments are presented in these papers:
a direct term-by-term
comparison of the infinite products defining the Selberg zeta; a
term-by-term comparison of the spectral expansions (via trace formula)
of the the logarithmic derivative of the Selberg zeta functions.
In generalizing the first approach to the $3$-dimensional case,
the obstacle consists in the more complicated structure of centralizers of hyperbolic
elements in discrete subgroups of $\SL(2,\C)$ than of $\SL(2,\R)$.
The difficulty in extending the spectral method to quotients of $\HH^3$
has to do with the continuous spectrum contribution to the trace formula.
For discussions of the Artin formalism in the context of the trace
formula, see \cite{Br-Sp},
and also, for the case of quotients of compact
manifolds, \cite{Jo-La}.

The proof of the identity of Selberg zeta functions
presented in this paper is very similar
to the one in \cite{Fr07}, which studied independently
the special case when $\Ga$ is a normal subgroup of $\Gn$.
The proof differs
from previous approaches in that it
avoids entirely the trace formula. Instead of comparing
traces of automorphic kernels,
we compare the (partial) automorphic kernels themselves.
This method can be extended easily
to quotients of real hyperbolic spaces of arbitrary dimension.

In the second part of the paper we prove the Artin formalism
for the scattering function by a direct method (independent of the Selberg
zeta function), and compare that to Venkov's original formula
in the 2-dimensional case.
\section{Definitions and notations}
Let $\HH^3:=\R^2\times(0,+\infty)$ be the upper half-space.
For an element $\om=(x_1,x_2,y)\in\HH^3$ we will
denote $y(\om)=y$ and $z(\om)=x_1+ix_2$.
$\HH^3$ carries $\SL(2,\C)$-invariant 
\newpage
\hspace*{-.4cm}metric
$ds^2=y^{-2}(dx_1^2+dx_2^2+dy^2)$ and hyperbolic volume
element $d\om=y^{-3}\,dx_1dx_2dy$.
The (positive) Laplace operator associated to this metric is
\vspace*{0.7cm}
\be
\Delta=-y^2\left(\frac{\dell^2}{\dell x_1^2}+\frac{\dell^2}{\dell x_2^2}
+\frac{\dell^2}{\dell y^2}\right)+y\frac{\dell}{\dell y}.
\ee
The \textbf{basic point-pair invariant} is given 
at \cite[Ch.1]{EGM} by
\be
\label{distance}
\delta(\om,\om')=\cosh\;d_{\HH^3}(\om,\om')=\frac{\|\om-\om'\|^2}{2y(\om) y(\om')}+1.
\ee
Assume $\Gamma\leq \SL(2,\C)$ is a discrete subgroup\footnote{We assume from
now on that all discrete subgroups of $\SL(2,\C)$ contain $\pm I_2$.}
such that the quotient $\Gamma\backslash\HH^3$ has finite hyperbolic volume.
Let $\chi$ be a finite dimensional unitary representation of $\Gamma$.
Then $\chi$ determines a vector bundle over $\GG\backslash\HH^3$, and we denote
by $L^2(\GG\backslash\HH^3,\chi)$ the Hilbert space of
$L^2$ sections.

Associated to
$\chi$, the Selberg zeta function $\ze$ is a function of
the complex parameter $s$
usually given
by an infinite product analogous to the Euler factorization of the Artin $L$-functions
(\cite{Se56, EGM, Fr05, Ve82}).
Its logarithmic derivative has a simpler expression and we will work with
that instead.

First, we need to introduce some notations.
(We follow \cite[Ch. 5]{EGM} and \cite{Fr05}.)
We denote the set of
hyperbolic\footnote{We will not need the distinction between hyperbolic and loxodromic elements in this paper.},
parabolic, and elliptic elements of $\GG$, respectively, by
\[\begin{split}
\Ghyp&:=\{P\in\GG: \tr(P)\notin [-2,2]\}
\\
\Gpar&:=\{P\in\GG: \tr(P)=\pm 2, P\neq\pm I_2\}
\\
\Gell&:=\{P\in\GG:\tr(P)\in (-2,2)\}
\end{split}\]
Note that we have a disjoint union
\[
\Ga=\{\pm I_2\}\cup\Ghyp\cup\Gpar\cup\Gell
\]
For $P\in\Ghyp$, $a(P)$ denotes the eigenvalue
with $|a(P)|>1$ and
$N(P)=|a(P)|^2$ is the norm of $P$.
Let $C_{\Ga}(P)$ denote the centralizer of $P$ in $\GG$.
This is a finitely generated abelian group of rank one. Let
\[C_{\Ga}(P)=C_{\Gamma}(P)^{\text{free}}\times
C_{\Ga}(P)^{\text{tor}}\] be
its decomposition into the free and torsion part.
A generator $P_0$ of $C_{\Gamma}(P)^{\text{free}}$
is called the primitive
hyperbolic element of $P$ in $\Gamma$.
Set $m(P):=|C_{\Ga}(P)^{\text{tor}}|$.

One can define uniquely the Selberg zeta function on $\ore s>1$ by
way of the logarithmic derivative
\be
\label{j1}
d\log\ze
=\sum_{[P]_{\GG}}\frac{\tr\chi(P)\log(NP_0)}{m(P)|a(P)-a(P)^{-1}|^2}
\,NP^{-s},
\quad \ore(s)>1
\ee
subject to the condition $\ddd\lim_{\ore s\to +\infty}Z_{\GG}(s,\chi)=1$.
The above sum is over $\GG$-conjugacy classes $[P]_{\GG}$ of hyperbolic
elements, with $P$ a representative of such a class. We use the notation
$\tr A$ for the trace of an operator $A$.
\begin{rek}
Let $\lam_j=1+t_j^2$, $j\geq 0$, the eigenvalues of
(the self-adjoint extension of)
$\Delta$ on
$L^2(\GG\backslash\HH^3, \chi)$.
With this definition, $\ze$ has zeros at $s=\pm it_j$, hence
on the line $\ore s=0$.
\end{rek}
\begin{defi}
Let $\ziu$ be the right-hand side of the
identity (\ref{j1}).
This is also the notation used in \cite{Ga-Wa}.
\end{defi}
It is known that $\ziu$
has a meromorphic continuation to
$s\in\C$ and its polar divisor is well understood, being related to the
spectrum of
$\Delta$ on $L^2(\GG\backslash\HH^3,\chi)$.
The meromorphic continuation
of $\ze$ itself is more problematic, in that it depends on
the specific shape of the polar divisor of $\ziu$.
Various instances are known when $\ze$
has meromorphic continuation to $s\in\C$: \cite{Ga-Wa}
in the general rank one case, when
$\Ga$ has no elliptic elements (and $\chi=1$);
\cite{Fr05} when $\Ga$ has no cuspidal elliptic elements (and arbitrary $\chi$).
Although, in most cases, \cite{Fr05} established the meromorphic
continuation of an integer power of $\ze$, the question is, in general, still
open (Conjecture in \cite[Section 6]{Fr05}).
\section{Main Results}
\subsection{The Artin formalism of the Selberg Zeta function}
Let $\Gn\subset\SL(2,\C)$ be a second Kleinian group containing $\Ga$
with finite index:
$\Ga\subset\Gn$ and $n=[\Gn:\Ga]$.
Let $\chi$ be a finite dimensional, unitary representation of $\Ga$, and
set $\pi=\Ind_{\Ga}^{\Gn}\chi$.
In this paper we prove
the following identity of meromorphic functions
\begin{thma}
\label{beef1}
For $s\in\C$,
\[\ziu=\ziun.\]
\end{thma}
As a corollary, we have
\be
\label{artin}
\ze=\zen,\quad\ore s>1,
\ee
and clearly this identity extends to the entire domain where the two Selberg zeta functions
admit meromorphic continuation.
\subsection{The Artin formalism of the scattering function}
Let $\scat$ denote the scattering function associated to
the continuous spectrum of $\Delta$ 
on $L^2(\Ga\backslash\HH^3,\chi)$.
(This is the determinant of the scattering matrix in the functional equation
of the Eisenstein series.)
It was observed earlier \cite{Ve-Zo} that the
identity (\ref{artin}) is linked, via the trace formula,
with a relationship
between $\scat$ and $\scatn$.
The second part of this paper is concerned with
a proof
of this relationship that can be developed separately
from (\ref{artin}),
by a direct matching of the Eisenstein series on the two spaces.
Namely, for a specific normalization
of the Eisenstein series, we have the following
\begin{thmb}
\label{beef2}
For $s\in\C$,
\[\scat=\scatn.\]
\end{thmb}
\section{Preliminaries}
\subsection{The map $\mapp$}
Let $\Cc(\Ga\backslash\HH^3)$ be the linear space of $\Ga$-invariant,
continuous functions
on $\HH^3$, and $\Bc(\Ga\backslash\HH^3)$ the subspace of
bounded, continuous functions. We use the analogous notation
for $\Gn$ instead of $\Ga$.

For $\{\alpha_i: 1\leq i\leq n\}$
a set of coset
representatives of $\Ga\backslash\Gn$
one can define the map
\be
\mapp:\Cc(\Ga\backslash\HH^3)\to \Cc(\Gn\backslash\HH^3),
\quad
\mapp f(\om)=\sum_{i=1}^nf(\alpha_i\, \om).
\ee
It can be easily checked that the map $\mapp$ is independent
of the particular choice of $\alpha_i$.
The following lemma is straightforward.
\begin{lemma}
\label{jint}
Let $f\in \Bc(\Ga\backslash\HH^3)$.
Then $\mapp f\in \Bc(\Gn\backslash\HH^3)$ and
\be
\int_{\Ga\backslash\HH^3}f(\om) d\om
=\int_{\Gn\backslash\HH^3}\mapp f(\om)d\om
\ee
\end{lemma}
\subsection{Orbital functions}
Assume $\Phi\in \Cc((1,+\infty))$ is a continuous function. This defines a 
point-pair invariant by
\be
k(\om,\om')\in \Cc(\HH^3\times\HH^3),\quad
k(\om,\om'):=\Phi\left(\delta(\om,\om')\right). 
\ee
Convergence issues set aside for the moment, by summing over all
the $\Ga$-translates of $k$ one obtains
the automorphic kernel
\be
K_{\chi}(\om,\om')=\sum_{\ga\in\Ga}k(\om,\gamma\om')\chi(\ga)
\ee
which is the kernel of an operator
$K_{\chi}$ on $L^2(\Ga\backslash\HH^3,\chi)$
(see \cite{Fr05} for details).
Instead of working with the automorphic kernel,
we will add translates over
proper subsets of $\GG$ to obtain $\GG$-automorphic functions, rather than kernels,
as follows.

Suppose that $\Om\subset \Ga$ is a
proper subset invariant under conjugation.
For our purposes we may assume that $\Om$ acts
without fixed points on $\HH^3$. That means
that $\Om$ does not contain
elliptic elements or the identity matrix.

Under the assumption that the function $\Phi$ has enough decay
at infinity, to each such a subset $\Om$ we can associate
the
$\Ga$-invariant function
\be
\label{orbital-definition}
F^{\Om}_{\chi}(\om):=\sum_{g\in\Om}k(\om,g\om)\tr\chi(g).
\ee
We will call such an object an orbital function,
as it is obtained by summing over disjoint unions of conjugacy orbits of $\Ga$.
Note that, formally, the orbital function associated to $\Om=\Ga$ is
the diagonal restriction of the automorphic kernel.

The next lemma is needed to establish the convergence of the
series defining $F_{\chi}^{\Om}$.
To state the lemma we first introduce some notation.
\begin{defi}
\label{cusps}
Let $h=h(\Ga)$ be the number of inequivalent cusps of $\Ga$,
and
$\germa_1,\dots, \germa_h$ a complete set of
$\Ga$-inequivalent cusps. For each $j$, $1\leq j\leq h$, we
fix $\si_j\in \SL(2,\C)$
such that $\si_j\cdot\infty=\germa_j$.
The $\GG$-invariant height function is defined as
\be
y_{\Ga}:\HH^3\to (0,+\infty),\quad
y_{\Ga}(\om):=\max_{1\leq j\leq h}
\max_{\gamma\in\GG}\;y(\si_j^{-1}\gamma\om).
\ee
\end{defi}
\noindent
The properties that make the height useful are:

i) $\ddd\inf_{\om\in\HH^3}y_{\Ga}(\om)>0$. 

ii) For a sequence $\{\om_n\}\subset\HH^3$,
$\{y_{\Ga}(\om_n)\}$ is unbounded
if and only if
$\{\om_n\}$ has a subsequence that converges to a cusp of $\GG$.
\begin{defi}
In general, if $\germa$ is an arbitrary cusp, we will use the notation
$\Ga_{\germa}$ for the stabilizer of $\germa$ in $\Ga$,
and $\Ga'_{\germa}:=\Ga_{\germa}\cap\Gpar$
for the subset of parabolic elements in the stabilizer.
\end{defi}
\begin{lemma}
\label{estimate}
Assume that $\si>0$ is a positive real number, and $\om\in\HH^3$.
Then, with $\delta(\om,\om')$ defined at (\ref{distance}), we have:
\asd
{\bf (a)}
$\ddd\sum_{\gamma\in\Gpar}
\delta(\gamma\om,\om)^{-\si}
\ll_{\si}\left[y_{\Ga}(\om)\right]^{2\si}$,\quad
for $\si>1$;
\asd
{\bf (b)}
$\ddd\sum_{\gamma\in\Ghyp}
\delta(\gamma\om,\om)^{-\si}
\ll_{\si}\left[y_{\GG}(\om)\right]^{2(2-\si)}$,\quad 
for $\si>2$.
\asd
Both estimates are uniform in $\om\in\HH^3$, with the implicit constant depending
on $\si$ (and $\GG$) only.
\end{lemma}
\emph{Proof.}
We assume for simplicity that $\Ga$ has only one cusp, namely at $\infty$.

Let $\Ginf$ be the stabilizer of $\infty$ in $\Ga$. Then
\be
\Ga'_{\infty}:=\Ginf\cap\Gpar
=\pm\left\{\mattwo{1}{\lam}{0}{1}: 0\neq\lam\in\Linf\right\},
\ee
where $\Linf$ is a lattice in $\C$ (see \cite[Thm. 2.1.8]{EGM}).

By Shimizu's lemma
(\cite[Prop. 2.3.7]{EGM}),
the  parabolic elements of $\Ga$
fix cusps , hence we have the disjoint union
\be
\label{Gpar}
\Gpar
=\bigcup_{g\in\Ga/\Ginf}\Ga'_{g\infty}
=\bigcup_{g\in\Ga/\Ginf}g\Ga'_{\infty}g^{-1}
=\bigcup_{g\in\Ginf\backslash\Ga}g^{-1}\Ga'_{\infty}g
\ee
Since $\delta\left(g^{-1}\mattwo{1}{\lam}{0}{1}g\,\om,\om\right)
=\delta(\lam+g\,\om,g\,\om)$
we have, for $\si>1$,
\[\begin{split}
&\sum_{\gamma\in\Gpar}\delta(\gamma\om,\om)^{-\si}
=\sum_{g\in\Ginf\backslash\Ga}\;
\sum_{0\neq\lam\in\Linf}\delta(\lam+g\om,g\om)^{-\si}
\\
&=\sum_{g\in\Ginf\backslash\GG}\;
\sum_{0\neq\lam\in\Linf}\left[\frac{|\lam|^2}{2y^2(g\om)}+1\right]^{-\si}
\ll \sum_{g\in\Ginf\backslash\GG}[y(g\om)]^{2\si}.
\end{split}\]
In the range $\si>1$
the latter series is,
by a standard estimate on the
Eisenstein series,
$O\left(y^{2\si}(\om)\right)$
as $y(\om)\to+\infty$.
This finishes the proof of part a) of the Lemma.

The proof of part b) follows the argument in \cite[p. 157]{EGM}.
Once again, we suppose for simplicity that $\GG$ has only one cusp.
Since $\Ghyp\subset\GG\setminus\Ginf$ we have, for $\si>2$,
\[\begin{split}
\sum_{\ga\in\Ghyp}\delta(\ga\om,\om)^{-\si}
&\leq\sum_{\ga\in\Ga-\Ginf}
\delta(\ga\om,\om)^{-\si}
=\sum_{t\in\Ginf}
\sum_{g\in\Ginf\backslash(\Ga-\Ginf)}
\delta(tg\om,\om)^{-\si}.
\end{split}\]
This is
\[\begin{split}
&\ll
\sum_{\lam\in\Linf}
\sum_{g\in\Ginf\backslash(\Ga-\Ginf)}
\delta(g\om+\lam,\om)^{-\si}
\\
&=
\sum_{g\in\Ginf\backslash(\Ga-\Ginf)}
\sum_{\lam\in\Linf}
\left[\frac{|\lam+z(g\om)-z(\om)|^2+y^2(\om)+y^2(g\om)}{2y(\om)y(g\om)}\right]^{-\si}
\\
&\ll
\sum_{g\in\Ginf\backslash(\Ga-\Ginf)}
y^{\si}(\om)\;y^{\si}(g\om)
\sum_{\lam\in\Linf}
\frac{1}{\;\left[\;|\lam+z(g\om)-z(\om)|^2+y^2(\om)\right]^{\si}}
\\
&\ll
\sum_{g\in\Ginf\backslash(\Ga-\Ginf)}
y^{\si}(\om)y^{\si}(g\om)
\cdot
[y^{-2\si}(\om)+y^{2-2\si}(\om)]\text{},
\end{split}\]
with the implied constant depending on $\si$ only.
As $y(\om)\to+\infty$, this is
\[
\ll y^{2-\si}(\om)
\sum_{g\in\Ginf\backslash(\Ga-\Ginf)}y^{\si}(g\om)\]
By the same standard estimate on the Eisenstein series, when $\si>2$,
\[\sum_{g\in\Ginf\backslash(\Ga-\Ginf)}y^{\si}(g\om)=O(y^{2-\si}(\om)),\quad
\text{as $y(\om)\to+\infty$},\]
which finishes the proof of part b).

{\it Note.}
In the general case (several cusps), one has to decompose the fundamental domain
into a finite union of cuspidal sectors (as in \cite[Thm. 1.2.4]{Ve82}) and 
analyze each cuspidal sector separately, but no new ideas are needed as far
as the estimates are concerned.
\begin{cor}
\label{bounds}
We fix the following assumptions:
\begin{itemize}
\item
$\Om\subset\Ga$ is a conjugation-invariant
subset, not containing elliptic elements or $\pm I_2$. That is,
$\Om\subset\Ghyp\cup\Gpar$;
\item
$\Phi\in \Cc ((1,\infty))$ is a continuous function 
with decay
$\Phi(x)=O(x^{-\alpha})$ as $x\to +\infty$.
\end{itemize}
Then:
\asd
{\bf (a)}
If $\alpha>1$, the series defining $F_{\chi}^{\Om}$ converges
absolutely and uniformly on compact subsets of $\HH^3$. Hence
$F_{\chi}^{\Om}\in \Cc(\Ga\backslash\HH^3)$;
\asd
{\bf (b)}
If $\alpha>2$ and $\Om\subset\Ghyp$,
then $F_{\chi}^{\Om}\in \Bc(\GG\backslash\HH^3)$.
\end{cor}

\begin{rek} The point here is that there is no requirement on
the behavior of $\Phi$ near $x=1$.
\end{rek}.

{\emph Proof.}
Fix a compact subset $K\subset\HH^3$.
Since $\Om$ does not fix any points in $\HH^3$,
\[\inf_{\ga\in\Om,\;\om\in K}\delta(\ga\om,\om)>1\]
Therefore, there exists a constant $C=C(K)$ such that
\[\Phi(\delta(\ga\om,\om))\leq C\delta(\ga\om,\om)^{-\alpha},\quad
\forall \ga\in\Om,\;\om\in K.\]
Using part a) of the Lemma we see now that the series defining
$F_{\chi}^{\Om}(\om)$ converges uniformly on compact subsets of $\HH^3$.

In the case $\Om\subset\Ghyp$, one can choose the constant $C$ uniformly
in $\om\in\HH^3$, hence for $\alpha>2$ we have
\[F_{\chi}^{\Om}(\om)\ll\sum_{\ga\in\Ghyp}\delta(\ga\om,\om)^{-\alpha}
\ll [y_{\Ga}(\om)]^{2(2-\alpha)}=O(1),\]
uniformly in $\om\in\HH^3$.
\subsection{}
Assume now that
$\Omn\subset\Gn$
is a conjugacy-invariant subset of $\Gn$, and
$\Om=\Omn\cap\Ga$.
With $\pi=\Ind_{\Ga}^{\Gn}$
we have the following
\begin{prop}
\label{orbital-prop}
Provided the series defining $F_{\chi}^{\Om}$ and $F_{\pi}^{\Omn}$
are uniformly convergent on compact sets, we have
\be
\label{orbital-identity}
\mapp F_{\chi}^{\Om}=F_{\pi}^{\Omn}.
\ee
\end{prop}
\emph{Proof.}
We follow the notation and the formula for the trace
of the induced representation from \cite[p.483]{Ve-Zo}.
Let
\[\tilde\chi(g)=\begin{cases}
\chi(g), &g\in\Ga\\
0, &g\notin\Ga.
\end{cases}\]
The trace of $\pi$ is then given by
\be
\tr\pi(g)=\sum_{i=1}^n\tr\tilde{\chi}(\alpha_i g\alpha_i^{-1}),\quad
g\in\Gn.
\ee
In the following computation we will be using two facts:
$\alpha_i\Omn\alpha_i^{-1}=\Omn$;
for $g\in\Omn$, $\tr\tilde\chi(g)=0$ unless
$g\in\Omn\cap\Ga=\Om$.
\[\begin{split}
F_{\pi}^{\Omn}(\om)
&=\sum_{i=1}^n\sum_{g\in\Omn}
k(\om,g\om)\cdot
\tr\tilde\chi(\alpha_i g\alpha_i^{-1})
\\
&=\sum_{i=1}^n\sum_{g\in\Omn}k(\om, \alpha_i^{-1}g\alpha_i\om)
\cdot \tr\tilde\chi(g)
=\sum_{i=1}^n\sum_{g\in\Om}k(\om,\alpha_i^{-1}g\alpha_i\om)
\cdot\tr\chi(g)
\\
&=\sum_{i=1}^n\sum_{g\in\Om}k(\alpha_i\om,g\alpha_i\om)
\cdot\tr\chi(g)
=\sum_{i=1}^nF_{\chi}^{\Om}(\alpha_i\om)
\\
&=\mapp F_{\chi}^{\Om}(\om)
\end{split}\]
By taking $f=F_{\chi}^{\Om}$ in Lemma \ref{jint} we obtain
\begin{cor}
\label{sdf}
Assume that the series defining $F_{\chi}^{\Om}$
is uniformly bounded
on $\HH^3$. Then
$F_{\chi}^{\Om}\in \Bc(\Ga\backslash\HH^3)$,
$F_{\pi}^{\Omn}\in \Bc(\Gn\backslash\HH^3)$,
and
\[\int_{\Ga\backslash\HH^3}F_{\chi}^{\Om}(\om)d\om
=\int_{\Gn\backslash\HH^3}F_{\pi}^{\Omn}(\om)d\om.\]
\end{cor}
In the next two sections,
we will show how the appropriate choice of $\Om$ and $\Phi$ translates
the identity of this Corollary into the result
stated at Theorem \ref{beef1}.
\section{Proof of Theorem 1.}
\subsection{Hyperbolic Orbital Functions}
In this section we restrict to the case $\Om=\Ghyp$,
and denote by
$F_{\chi}^{\rm hyp}:=F_{\chi}^{\Ghyp}$ the corresponding orbital function.
By Proposition \ref{bounds},
$F_{\chi}^{\rm hyp}\in \Bc(\Ga\backslash\HH^3)$
whenever $\Phi(x)=O(x^{-2-\ep})$, for some $\ep>0$.
The computation of its integral is standard
in the theory of Selberg trace formula
(see \cite[Ch.5]{EGM}):
\be
\label{hyperbolic}
\begin{split}
\int_{\Ga\backslash\HH^3}F_{\chi}^{\rm hyp}(\om)d\om
&=\sum_{[P]_{\GG}}\tr\chi(P)
\int_{\Ga_P\backslash\HH^3}k(\om, P\om)\,d\om
\\
&=\sum_{[P]_{\GG}}\frac{\tr\chi(P)\; \log(NP_0)}{m(P)|a(P)-a(P)^{-1}|^2}
\; g(\log NP).
\end{split}
\ee
Here the sum is over hyperbolic $\Ga$-conjugacy classes
$[P]_{\GG}$, $\Ga_P$ is the stabilizer of $P$ in $\Ga$,
and $g(t)$ is
obtained from $\Phi(x)$ via the Selberg transform
\be
\label{transform}
g(t)=2\pi\int_{\cosh t}^{\infty}\Phi(x)\,dx.
\ee
\begin{rek}
Under more restrictive conditions on the test function $\Phi(x)$,
the right-hand side of the identity (\ref{hyperbolic}) represents
precisely the contribution of the hyperbolic conjugacy classes
to the geometric side of the trace formula.
\end{rek}
\subsection{Green function kernel}
For $s\in\C$, with $\ore(s)>1$, we consider, as in
\cite[p.185]{EGM},
the test function
\be
\Phi_s(x)=\frac{2^{-s}s}{\pi}\cdot
\frac{(x+\sqrt{x^2-1})^{-s}}{\sqrt{x^2-1}}.
\ee
This is a smooth function on $(1,\infty)$, and clearly
$\Phi_s(x)=O\left(x^{-\ore(s)-1}\right)$ as $x\to +\infty$.
The corresponding transform given by (\ref{transform}) is
\be
\label{g-formula}
g_s(t)=e^{-s|t|}.
\ee
We fix now the parameter $s\in\C$, with $\ore s>1$.
\noindent
Let $F_{\chi}^{\rm hyp}(\om;s)$ be
the corresponding orbital function on $\GG\backslash\HH^3$
associated to $\chi, \Ghyp, \Phi_s$.
Hence
$F_{\chi}^{\rm hyp}(\om;s)\in \Bc(\Ga\backslash\HH^3)$
and the formula (\ref{hyperbolic}) yields
\be
\label{hyp2}
\int_{\Ga\backslash\HH^3}
F_{\chi}^{\rm hyp}(\om;s)d\om=\ziu,\quad
\ore(s)>1.
\ee
where $\ziu$ has been defined at (\ref{j1}).
\subsection{Conclusion}
With the choice of test function $\Phi=\Phi_s$, $\ore s>1$,
and the choice of conjugacy-invariant subsets
$\Om=\Gpar\subset\Ga$ and $\Omn=\Gn_{\rm par}\subset\Gn$,
the result of Cor. \ref{sdf},
\[\int_{\Ga\backslash\HH^3}F_{\chi}^{\rm hyp}(\om;s)d\om
=\int_{\Gn\backslash\HH^3}F_{\pi}^{\rm hyp}(\om;s)d\om,\]
combined with the identity (\ref{hyp2}), yields
\be
\ziu=\ziun, \quad\ore s>1.
\ee
\section{Proof of Theorem 2.}
It was observed in \cite{Ve-Zo} that the identity
stated at Theorem 2 can be obtained as a by-product of the
Artin formalism for the Selberg zeta functions combined with the Selberg trace
formula.

We present here a direct proof which is independent of the identity
of Selberg zeta functions.
For the simplicity of exposition we restrict ourselves, throught
this section only, to the case when
\be
\label{assumption}
\textit{$\Gn$ has only one cusp (at $\infty$) and $\dim\chi=1$,}
\ee
but we remark that our argument works in complete generality
(see Remark \ref{rek:generalization} below).

For $1\leq i\leq h=h(\Ga)$, let $\Ga_i:=\Ga_{\germa_i}$ the stabilizer of
the cusp $\germa_i$ in $\Ga$, and $\Gn_i=\Gn_{\germa_i}$
the stabilizer of $\germa_i$ in $\Gn$.
We set
\be
n_i=[\Gn_i:\Ga_i],\quad 1\leq i\leq h
\ee
Since all the cusps of $\Gn$ are equivalent to $\infty$,
we will choose the scaling matrices (introduced in Definition \ref{cusps})
such that $\si_i\in\Gn$.
\subsection{Eisenstein series associated to $\chi$}
The degree of singularity of $\chi$
(cf. \cite{Ve-Zo}, \cite{Fr05})
is the integer $1\leq\kappa\leq h$
with the following property:
the restriction of $\chi$ to $\Ga_j$ is trivial for
$1\leq j\leq\kappa$, and non-trivial for $\kappa+1\leq j\leq h$.

The Eisenstein series on $\Ga\backslash\HH^3$
associated to $\chi$ are defined, for $\ore s>2$,
by the absolutely convergent series
\be
E_j(\om,s;\chi)
:=\sum_{g\in\Ga_j\backslash\Ga}
y^s(\si_j^{-1}g\cdot\om)\chi(g^{-1})
\in\mathscr A(\Ga,\chi),\quad 1\leq j\leq\kappa.
\ee
Here $\mathscr A(\Ga,\chi)$ is the linear space of $\chi$-automorphic forms
of polynomial growth on $\GG\backslash\HH^3$
(see \cite[Chap. 3]{EGM} and \cite{Fr05}, for a more precise definition).

It is known that the Eisenstein series have
meromorphic continuation to $s\in\C$ (\cite{Co-Sa, EGM, Fr05}).

Given a choice
of coset representatives $\{\alpha_i: 1\leq i\leq n\}$ of
$\Ga\backslash\Gn$,
it was shown in \cite{Ve-Zo} that the representation
$\pi=\Ind_{\Ga}^{\Gn}\chi$ can be realized on
$V:=\C^n$ (column vectors) via left multiplication by the $n\times n$ matrices
\be
\pi(g)=\left[\tilde\chi\left(\alpha_ig\alpha_j^{-1}\right)\right]_{1\leq i, j \leq n},
\quad g\in\Gn.
\ee
\subsection{Coset decomposition}
We will use a specific choice of representatives described by the following
proposition.
\begin{prop}
For $1\leq i\leq h$, let
$\{\beta_{it}: 1\leq t\leq n_i\}$ be
a set of coset representatives of $\Ga_i\backslash\Gn_i$.
Then $\Gn=\bigcup_{i=1}^h\bigcup_{t=1}^{n_i}
\Ga_i\beta_{it}\si_i$ is
a disjoint union.
That is, we can take
$\{\alpha_\nu: 1\leq \nu\leq n\}
=\{\beta_{it}\si_i: 1\leq i\leq h, 1\leq t\leq n_i\}$.
In particular,
\[\sum_{i=1}^hn_i=n.\]
\end{prop}
\emph{Proof.}
Let $g'\in\Gn$. There exists $i$, $1\leq i\leq h$, and $g\in\Ga$
such that $g'\infty=g\germa_i$.
Since $\infty=\si_i^{-1}\germa_i
\Rightarrow g'\si_i^{-1}\cdot\germa_i=g\germa_i
\Rightarrow g^{-1}g'\si_i^{-1}\in\Gn_i$.
Consequently $g'\in\Ga\Gn_i\si_i$.
But $\Gn_i=\bigcup_t\Ga_i\beta_{it}$,
hence
\be
g'\in\bigcup_{t=1}^{n_i}\Ga\Ga_i\beta_{it}\si_i
\subset\bigcup_{t=1}^{n_i}\Ga\beta_{it}\si_i
\subset
\bigcup_{i=1}^{h}\bigcup_{t=1}^{n_i}
\Ga\beta_{it}\si_i
\ee
It is then straightforward to check that the cosets
$\Ga\beta_{it}\si_i$
are actually disjoint for different pairs of $(i,t)$.

From now on we will identify the sets
$\{\alpha_{\nu}\}$ and
$\{\beta_{it}\si_i\}$
in the following order:
\be
\label{indices}
\alpha_{\nu}=\beta_{it}\si_i,
\quad{\rm if} \quad \nu=n_1+\dots+n_{i-1}+t
\ee
\subsection{Eisenstein series associated to $\pi$}
The singular space of the induced representation $\pi$ is
\be
V_{\infty}:=\{v\in V: \pi(g)v=v, \quad\forall g\in\Ginf\}.
\ee
It is known that it has dimension $\kappa$ and orthonormal basis
\be
\eje=\frac{1}{\sqrt n_j}[0,\dots 0, 1\dots, 1,0,\dots,0]^t,\quad
1\leq j\leq \kappa,
\ee
where the $1$'s occur in the $j^{\rm th}$ block (of length $n_j$), according to the 
identification (\ref{indices}).
The Eisenstein series on $\Gn\backslash\HH^3$ associated to $\pi$ are
given, for $\ore s>2$, by
\be
E_j(\om,s;\pi)
=\sum_{g\in\Gn_{\infty}\backslash\Gn}y^s(g\cdot \om)
\pi(g^{-1})\eje\in
\mathscr A(\Gn,\pi),\quad 1\leq j\leq\kappa.
\ee
\subsection{The map $\mapt$}
In \cite{Ve-Zo} the following map is defined:
\be
\mapt:\mathscr A(\Ga,\chi)\to\mathscr A(\Gn,\pi),\quad
\mapt f(\om)=[f(\alpha_{1}\om),\dots, f(\alpha_n\om)]^t.
\ee
The following proposition represents the main result of this section.
\begin{prop}
\label{eiscor}
For $1\leq j\leq \kappa$,
\[\mapt E_j(\om,s;\chi)=n_j^{1/2}\,E_j(\om,s;\pi).\]
\end{prop}
\emph{Proof.}
Let $\nu=n_1+\dots+n_{i-1}+t$,
that is, $\alpha_{\nu}=\beta_{it}\si_i$.
The $\nu^{\rm th}$ component of $E_j(\om,s;\pi)$,
as a vector in $\C^n$, is
\[
\frac{1}{\sqrt n_j}
\sum_{g\in\Gn_{\infty}\backslash\Gn}
y^s(g\cdot\om)\sum_{k=1}^{n_j}
\tilde\chi
\left(\beta_{it}\si_ig^{-1}(\beta_{jk}\si_j)^{-1}\right).
\]
With the change of variable $g\mapsto\si_jg$ and using
the fact that $\si_j\Gn_{\infty}\si_j^{-1}=\Gn_j$, the above sum becomes
\[
\frac{1}{\sqrt n_j}\sum_{g\in\Gn_j\backslash\Gn}
y^s(\si_j^{-1}g\cdot\om)\sum_{k=1}^{n_j}
\tilde\chi\left(\beta_{it}\si_ig^{-1}\beta_{jk}^{-1}\right).
\]
With a further change of variable
$g\mapsto \beta_{jk}g$, this equals:
\[\begin{split}
&\frac{1}{\sqrt n_j}
\sum_{g\in\Gn_j\backslash\Gn}
y^s(\si_j^{-1}\beta_{jk}^{-1}g\cdot\om)
\sum_{k=1}^{n_j}
\tilde\chi\left(\beta_{it}\si_ig^{-1}\right)
\\
&=\frac{1}{n_j^{3/2}}
\sum_{g\in\Ga_j\backslash\Gn}
\sum_{k=1}^{n_j}\;
\dots\qquad\text{[summing over $\Ga_j\backslash\Gn$ instead of
$\Gn_j\backslash\Gn$]}
\\
&=\frac{1}{n_j^{3/2}}
\sum_{g\in\Ga_j\backslash\Ga}
\sum_{\theta\in\Ga\backslash\Gn}
\sum_{k=1}^{n_j}
y^s(\si_j^{-1}\beta_{jk}^{-1}g\theta\cdot\om)\,
\tilde\chi
\left(\beta_{it}\si_i\theta^{-1}\right)\chi(g^{-1})
\end{split}\]
Now,
$\tilde\chi(\beta_{it}\si_i\theta^{-1})=0$ unless
$\beta_{it}\si_i\theta^{-1}\in\Ga$. This means that
$\Ga\beta_{it}\si_i=\Ga\theta$, which
forces $\theta$ in the sum
$\sum_{\theta\in\Ga\backslash\Gn}$ to equal
$\beta_{it}\si_i$.
Hence the $\nu^{\rm th}$ component of $E_j(\om,s;\pi)$ equals
\[\begin{split}
&\frac{1}{n_j^{3/2}}
\sum_{g\in\Ga_j\backslash\Ga}
\sum_{k=1}^{n_j}
y^s(\si_j^{-1}\beta_{jk}^{-1}g\beta_{it}\si_i\cdot \om)\,
\chi(g^{-1}).
\end{split}\]
Note that
\[y(\si_j^{-1}\beta_{jk}^{-1}\cdot P)
=y(\si_j^{-1}\beta_{jk}^{-1}\si_j\cdot\si_j^{-1}P)
=y(\si_j^{-1}P),\quad P\in\HH^3,\]
since $\si_j^{-1}\beta_{jk}^{-1}\si_j\in\Gn_{\infty}$.
Therefore the terms in the $k$-sum are all equal, and
the $\nu^{\rm th}$ component of $E_j(\om,s;\pi)$ equals
\[\begin{split}
&\frac{1}{n_j^{3/2}}
\sum_{g\in\Ga_j\backslash\Ga}
\sum_{k=1}^{n_j}
y^s(\si_j^{-1}g\beta_{it}
\si_i\cdot\om)\,\chi(g^{-1})
\\
&=\frac{1}{\sqrt n_j}
\sum_{g\in\Ga_j\backslash\Ga}
y^s(\si_j^{-1}g\beta_{it}
\si_i\cdot\om)\,\chi(g^{-1})
\\
&=\frac{1}{\sqrt n_j}
E_j(\beta_{it}\si_i\cdot\om,s;\chi)
=\frac{1}{\sqrt n_j}
E_j(\alpha_{\nu}\om,s;\chi) \qquad
\text{[recall that $\alpha_{\nu}=\eta_{it}\si_i$]}.
\end{split}\]
This shows that $E_j(\om,s;\pi)=\frac{1}{\sqrt n_j}\mapt E_j(\om,s;\chi)$.
\begin{rek}\label{rek:generalization}
The general version of Proposition \ref{eiscor}
(for an arbitrary number
of cusps of $\Gn$ and $\chi$ of arbitrary dimension)
was stated and proved in \cite[Prop. 14]{Br-Sp}.
The proof of Theorem 2, which is an immediate corollary
of Prop. \ref{eiscor}
(as shown in the next section),
thus carries over to the general case.
\end{rek}
\subsection{The Artin formalism of the scattering function}
Let
\[\mathcal E(\om,s;\chi)=[E_1(\om,s;\chi),\dots, E_{\kappa}(\om,s;\chi)]^t\]
be the column vector that encodes all the Eisenstein series
associated to $\chi$. The scattering matrix
$\mathfrak S_{\Ga}(s;\Ga;\chi)
=\left[\mathfrak S_{ij}(s;\chi)\right]_{1\leq i,j\leq \kappa}$ is determined by the functional equation (we refer to \cite{Fr05} for more details)
\[
\mathcal E(\om,s;\chi)=\mathfrak S(s;\Ga;\chi)
\mathcal E(\om,2-s;\chi).\]
That is,
\be
\label{scatter}
E_i(\om,s;\chi)=\sum_{j=1}^{\kappa}
\mathfrak S_{ij}(s;\chi)E_j(\om,2-s;\chi),\quad
1\leq i\leq\kappa.
\ee
The scattering function is
\be
\scat:=\det\mathfrak S(s;\Ga;\chi).
\ee

Similarly, the scattering matrix of $\pi$ is determined by
\be
\label{scatter2}
E_i(\om,s;\pi)=\sum_{j=1}^{\kappa}
\mathfrak S_{ij}(s;\pi)E_j(\om,2-s;\pi),\quad
1\leq i\leq\kappa.
\ee
By applying the $\mapt$ operator on both sides of the equation (\ref{scatter})
we obtain, in view of Proposition \ref{eiscor},
\[
n_i^{1/2}E_i(\om,s;\pi)
=\sum_{j=1}^{\kappa}
\mathfrak S_{ij}(s;\chi)n_j^{1/2}E_j(\om,s;\pi).
\]
Comaparison to (\ref{scatter2}) gives the relation
between the two scattering matrices:
\be
\mathfrak S_{ij}(s;\pi)=n_i^{-1/2}\mathfrak S_{ij}(s;\chi)n_j^{1/2},
\quad
\text{hence}
\quad
\mathfrak S(s;\Gn;\pi)=D^{-1}\mathfrak S(s;\Ga;\chi)D,
\ee
where $D$ is the $\kappa\times\kappa$ diagonal matrix
with $n_i^{1/2}$ on the diagonal. This implies
that the two scattering matrices have the same
determinant, which is the statement of Theorem 2.
\subsection{Comparison with Thm. 3.2 in \cite{Ve-Zo}}
The aim of this section is to reconcile Theorem 2 with
the slightly different version of Artin formalism proved in \cite{Ve-Zo}
for the 2-dimensional case.

In the remaining section
$\Ga\subset\Gn\subset \SL(2,\R)$
are cofinite Fuchsian groups.
The result (and proof) of Theorem 2 carries
over to this case almost word by word:
if $\chi$ is a finite dimensional
unitary representation and
$\pi=\Ind_{\Ga}^{\Gn}\chi$, then
$\scat=\scatn$.
However, Theorem 3.2 in \cite{Ve-Zo} states the following
(we will use the upper index $VZ$ for the analogous concepts introduced in $VZ$):
\be
\label{VZ}
\phi^{VZ}(s;\Ga;\chi)\Om(\chi)^{1-2s}=\phi^{VZ}(s;\Gn;\pi)\Om(\pi)^{1-2s},
\ee
(with the constants $\Om(\pi)$ and $\Om(\chi)$ to be defined shortly).
The goal of this section is a direct proof
that the two formulas are equivalent.

We will keep the assumption that $\Gn$ has only one cusp
(at $\infty$) and $\dim\chi=1$, for simplicity of exposition.
We keep the notation $\germa_1,\dots, \germa_h$ for the cusps of
$\Ga$.

In the 2-dimensional case, stabilizers of cusps in Fuchsian groups are cyclic
($\mmod \pm I_2$). Let $\Sn$ be a generator of $\Gn_{\infty}$. It is easy
to see that $S_i:=\si_i\Sn^{n_i}\si_i^{-1}$ is a generator
of $\Ga_i:=\Ga_{\germa_i}$, for $1\leq i\leq h$.

The constants $\Om(\chi)$ and $\Om(\pi)$ are defined as follows. First,

\be
\Om(\chi):=\prod_{j=1}^{\kappa}|1-\chi(S_j)|.
\ee
Let now
$\det'(A)$ denote the product of the non-zero eigenvalues of a
matrix $A$. Then
\be
\Om(\pi):=\big|\det{}'(I-\pi(\Sn))\big|,
\ee
with $I$ the identity matrix.

The difference between the two versions of the Artin formalism
for the scattering function
stems from different choices of
of the scaling matrices of cusps (Def. \ref{cusps}).
This in turn leads to a different
normalization of the Eisenstein series (associated to $\chi$),
and thus to an extra scalar factor
in the scattering function. In \cite{Ve-Zo} the scaling matrices
of cusps (see Def. \ref{cusps}) are chosen such that
the lattices associated to stabilizers of cusps
have (relative) co-volume one. That is,
\[\si_j^{VZ}=
\si_j\cdot\mattwo{n_j^{1/2}}{}{}{n_j^{-1/2}},\]
which leads to
\[E_j^{VZ}(\om,s;\chi)
=n_j^{-s}E_j(\om,s;\chi),
\quad 1\leq j\leq \kappa.\]
Therefore, with $D$ the $\kappa\times\kappa$ diagonal matrix with $n_j$ on the diagonal, we
have
\[\begin{split}
&\mathcal E^{VZ}(\om,s;\chi)
=D^{-s}\mathcal E(\om,s;\chi)
\\
&=D^{-s}\mathfrak S(s;\Ga;\chi)\mathcal E(\om,1-s;\chi)
\qquad\text{[functional equation in the 2-dimensional case]}
\\
&=D^{-s}\mathfrak S(s;\Ga;\chi)D^{1-s}
\mathcal E^{VZ}(\om,s;\chi),
\end{split}\]
which gives
\be
\mathfrak S^{VZ}(s;\Ga;\chi)
=D^{-s}\mathfrak S(s;\Ga;\chi)D^{1-s}.
\ee
The scattering determinants then satisfy the relation
\be
\phi_{\Ga}^{VZ}(s;\Ga;\chi)=
\scat\cdot
\prod_{j=1}^{\kappa}n_j^{1-2s}
\ee
On the other hand, $E_j^{VZ}(\om,s;\pi)=E_j(\om,s;\pi)$, hence
$\phi^{VZ}(s;\Gn;\pi)=\scatn$.

Therefore, to show the equivalence of the two
versions of Artin formalism for the scattering function,
we will prove the following formula.
\begin{prop}
\label{constants}
With the above notations,
\[\frac{\Om(\pi)}{\Om(\chi)}=\prod_{j=1}^{\kappa}n_j.\]
\end{prop}
\emph{Proof.} As remarked before, in the 2-dimensional case stabilizers of cusps
are cyclic, hence we have the following coset representatives for
$\Ga_i\backslash\Gn_i$:
\be
\beta_{ia}=\si_i\Sn^a\si_i^{-1},\quad 1\leq a\leq n_i,\quad 1\leq i\leq h.
\ee
Recall that
\[\pi(\Sn)=[\tilde\chi(\alpha_i\Sn\alpha_j^{-1}]_{1\leq i,j\leq n}\]
where $\{\alpha_{\nu}:1\leq\nu\leq n\}
=\{\beta_{ia}\si_i: 1\leq i\leq h, 1\leq a\leq n_i\}$.
It is easy to verify that
$\beta_{ia}\si_i\Sn(\beta_{jb}\si_j)^{-1}\notin\Ga$, unless $i=j$.
Since $\beta_{ia}\si_i\Sn(\beta_{ib}\si_i^{-1})^{-1}=\si_i\Sn^{a+1-b}\si_i^{-1}$,
we find that
$\pi(\Sn)$ is a block-diagonal matrix
\be
\pi(\Sn)=\bigoplus_{i=1}^{h}\;
[\,\tilde\chi(\si_i\,\Sn^{a+1-b}\,\si_i^{-1})\,]_{1\leq a,b\leq n_i}
\ee
Moreover, $\si_i\Sn^{a+b-1}\si_i^{-1}\in\GG$
if and only if $a+b-1\equiv 0\,(\text{mod}\, n_i)$. Hence
\be
\pi(\Sn)=\bigoplus_{j=1}^h
\left(\begin{array}{ccccc}
0&1&0&\dots&0\\
0&0&1&\dots&0\\
\dots&\dots&\dots&\dots&\dots\\
0&0&\dots&0&1\\
\chi(S_j)&0&\dots&0&0
\end{array}\right)
\ee
A straightforward computations gives the characteristic polynomial of
$\pi(\Sn)$:
\be
P(\lam)=\prod_{j=1}^h(\lam^{n_j}-\chi(S_j))
\ee
Therefore the eigenvalue $\lam=1$ comes only from the first $\kappa$ blocks,
where $\chi(S_j)=1$.
The product of non-zero eigenvalues of $I-\pi(\Sn)$ is
\[\begin{split}
\det{}'(I-\pi(\Sn))
&=\prod_{1\leq j\leq \kappa}
\left[\frac{\lam^{n_j}-1}{\lam-1}\Big|_{\lam=1}\right]
\cdot\prod_{i=\kappa+1}^{h}(1-\chi(S_i))
\\
&=\prod_{j=1}^{\kappa}n_j\cdot\prod_{i=\kappa+1}^h(1-\chi(S_j))
\end{split}\]
Taking the absolute value on both sides of this identity
finishes the proof of Proposition \ref{constants}.

\vspace*{.5cm}
\emph{Acknowledgments.}  We are grateful to Professor Jay Jorgenson for helpful discussions.  The first-named author thanks the Center for the Advanced Study of Mathematics at Ben-Gurion University for supporting him while the research leading to this paper was carried out.  He also thanks Mr. Tony Petrello for additional travel
support that made the collaboration possible.

\end{document}